\documentclass[11pt]{article}
\usepackage{amsmath}
\usepackage{amssymb}
\usepackage{amsfonts}

\makeatletter

\@addtoreset{equation}{section}
\makeatother
\def\<{\langle}
\def\>{\rangle}

\newtheorem{lem}{Lemma}[section]
\newtheorem{theo}{Theorem}[section]
\newtheorem{rem}{Remark}[section]

\makeatletter
   
   \@addtoreset{equation}{section}
\makeatother

\begin{document}
\title{\bf Asymptotic Profile of Solutions for Strongly Damped Klein-Gordon Equations}
\author{Ryo IKEHATA\thanks{Corresponding author: ikehatar@hiroshima-u.ac.jp} \\ {\small Department of Mathematics, Graduate School of Education, Hiroshima University} \\ {\small Higashi-Hiroshima 739-8524, Japan}}
\maketitle
\begin{abstract}
We consider the Cauchy problem in ${\bf R}^{n}$ for strongly damped Klein-Gordon equations. We derive asymptotic profiles of solutions with weighted $L^{1,1}({\bf R}^{n})$ initial data by a simple method introduced by R. Ikehata. The obtained results show that the wave effect will be weak because of the mass term, especially in the low dimensional case ($n = 1,2$) as compared with the strongly damped wave equations without mass term ($m = 0$), so the most interesting topic in this paper is the $n = 1,2$ cases.
\end{abstract}

\section{Introduction.}
\footnote[0]{Keywords and Phrases: Klein-Gordon equation; Structural damping; Asymptotic profiles; Fourier Analysis; Low frequency; High frequency; Weighted $L^{1}$-initial data; Low dimension.}
\footnote[0]{2010 Mathematics Subject Classification. Primary 35L15, 35L05; Secondary 35B40, 35B65.}
We are concerned with the Cauchy problem for the so called Klein-Gordon type of equation in ${\bf R}^{n}$ ($n \geq 1$) with the structural damping and mass terms
\begin{equation}
u_{tt}(t,x) - \Delta u(t,x)  + m^{2}u(t,x) -\Delta u_{t}(t,x) = 0,\ \ \ (t,x)\in (0,\infty)\times {\bf R}^{n} ,\label{eqn}
\end{equation}
\begin{equation}
u(0,x)= u_{0}(x), \quad u_{t}(0,x) = u_{1}(x), \quad x\in {\bf R}^{n},\label{initial}
\end{equation}
where $m > 0$. The initial data $u_{0}$ and $u_{1}$ are also chosen from the usual energy space
\[[u_{0},u_{1}] \in H^{1}({\bf R}^{n}) \times L^{2}({\bf R}^{n}).\]
Then, it is well-known (cf. see \cite{ITY}) that the problem (1.1)-(1.2) has a unique weak solution
\[u \in C([0,\infty);H^{1}({\bf R}^{n})) \cap C^{1}([0,\infty);L^{2}({\bf R}^{n})).\]
\noindent
We first mention several known facts concerning the strongly damped wave equations without mass term (i.e., $m = 0$) such that
\begin{equation}
u_{tt}(t,x) - \Delta u(t,x)  -\Delta u_{t}(t,x) = 0,\ \ \ (t,x)\in (0,\infty)\times {\bf R}^{n}.
\end{equation}
At the first stage, Ponce \cite{p} and Shibata \cite{Shibata} studied the decay structure of the $L^{p}$-norm of solutions (including their derivatives) to the problem (1.3) and (1.2). Recently, by the results from Ikehata-Todorova-Yordanov \cite{ITY}, Ikehata \cite{Ik-4} and Ikehata-Onodera \cite{IO} one can know that the asymptotic profile of the solution to the equation (1.3) is the so called diffusion waves, that is,
\[\hat{u}(t,\xi) \sim P_{1}e^{-t\vert\xi\vert^{2}/2}\frac{\sin(t\vert\xi\vert)}{\vert\xi\vert} \quad (t \to \infty),\]
where 
\[P_{j} := \displaystyle{\int_{{\bf R}^{n}}}u_{j}(x)dx,\quad (j = 0,1),\]
and 
\[\hat{u}(t,\xi) := {\cal F}(u(t,\cdot))(\xi) := \frac{1}{\sqrt{(2\pi)^{n}}}\int_{{\bf R}^{n}}e^{-ix\cdot\xi}u(t,x)dx\]
represents the partial Fourier transform of the solution $u(t,x)$ with respect to the $x$-variable. In this connection, concerning the higher order asymptotic expansion of the solution to problem (1.3) and (1.2), one can cite quite recent results due to Michihisa \cite{M} (see also Ikehata-Takeda \cite{IT} for an  asymptotic behavior of solutions for a more generalized equations of (1.3) from a different point of view). Furthermore, as an application to the nonlinear problem of the equation (1.3) one can cite D'Abbicco-Reissig \cite{DR} and the references therein. Anyway, from the results obtained by \cite{Ik-4} and \cite{IO} one can know that as $t \to \infty$, in the case when $n \geq 3$ one has
\[C\vert P_{1}\vert t^{-\frac{n-2}{4}} \leq \Vert u(t,\cdot)\Vert \leq C^{-1}(\Vert u_{1}\Vert_{1,1} + \Vert u_{0}\Vert_{1,1} + \Vert u_{1}\Vert + \Vert u_{0}\Vert_{H^{1}})t^{-\frac{n-2}{4}},\]
while in the case when $n = 1,2$ it is true that
\[C\vert P_{1}\vert \sqrt{\log t} \leq \Vert u(t,\cdot)\Vert \leq C^{-1}(\Vert u_{1}\Vert_{1,1} + \Vert u_{0}\Vert_{1,1} + \Vert u_{1}\Vert + \Vert u_{0}\Vert_{H^{1}})\sqrt{\log t}, \quad (n = 2)\]
\[C\vert P_{1}\vert \sqrt{t} \leq \Vert u(t,\cdot)\Vert \leq C^{-1}(\Vert u_{1}\Vert_{1,1} + \Vert u_{0}\Vert_{1,1} + \Vert u_{1}\Vert + \Vert u_{0}\Vert_{H^{1}})\sqrt{t}, \quad (n = 1).\]
Interestingly, through the solution to the strongly damped wave equation (1.3), (at least) in the one and two dimensional whole space cases we can find that the so called Poincar\'e inequality never hold, because one can also get the decay estimates of the total energy (see Ikehata-Natsume \cite[Theorem 1.1]{INatsume}, and also Char\~ao-da Luz-Ikehata \cite{CLI}) such that for all $n \geq 1$
\[\Vert u_{t}(t,\cdot)\Vert + \Vert\nabla u(t,\cdot)\Vert \leq C(1+t)^{-\frac{n}{4}}\Vert u_{1}\Vert_{1} + C(1+t)^{-\frac{n+2}{4}}\Vert u_{0}\Vert_{1} + Ce^{-\alpha t}(\Vert u_{1}\Vert + \Vert\nabla u_{0}\Vert).\]
This phenomena is coming from an observation on the asymptotic profile itself such that in the low dimensional case $n = 1,2$ the effect of the wave part of the profile coming from the factor $\sin(t\vert\xi\vert)/\vert\xi\vert$ is extremely strong, so that the diffusive structure of the solution vanishes any more as time goes to infinity. This interesting phenomena shows a quite different aspect as compared with the usual weakly damped wave equation case, which is nowadays very popular due to many papers published by D'Abbicco-Ebert \cite{DE}, Hosono-Ogawa \cite{HO}, Karch \cite{K}, Narazaki \cite{Na}, Nishihara \cite{N-2}, Takeda \cite{T}, Wakasugi \cite{W} and the references therein:
\begin{equation}
u_{tt}(t,x) - \Delta u(t,x) + u_{t}(t,x) = 0,\ \ \ (t,x)\in (0,\infty)\times {\bf R}^{n}.
\end{equation} 
In the weakly damped wave equation case (1.4) the $L^{2}$-norm of solutions necessarily decays when the total energy decays with some rate. This is because of the diffusive aspect of the equation (1.4). In this sense, the strongly damped wave equation (1.3) does not have a diffusive structure anymore at least in the low dimensional case $n = 1,2$.

In this connection, from these observations one can find that when we want to get the "decay" property of the solution to the Cauchy problem of the equation (1.3), we have to impose rather a stronger assumption on the initial velocity $u_{1}(x)$ such that $P_{1} = 0$. People sometimes impose this zero average condition on the initial data when they study the decay property of solutions to the equation (1.4) in the case when $n = 2$ (and/or $n = 1$). However, this zero average condition $P_{j} = 0$ ($j = 0,1$) seems to be just a technical one (at least) in the case of (1.4), while in the case of (1.3) such a zero average condition $P_{1} = 0$ seems to be essential in the low dimensional case $n = 1,2$ in order to get the $L^{2}$-decay of solutions. It should be emphasized that the paper due to D'Abbicco-Ebert-Picon \cite{DEP} have also pointed out its importance of the zero average condition $P_{j} = 0$ ($j = 0,1$) in order to get the decay estimates of solutions to the equations 
\[u_{tt}(t,x) + (- \Delta)^{\sigma}u(t,x) + (- \Delta)^{\delta}u_{t}(t,x) = 0,\]
in terms of Hardy spaces.\\

If one considers the following weakly damped Klein-Gordon equation in place of (1.1):
\begin{equation}
u_{tt}(t,x) - \Delta u(t,x)  + m^{2}u(t,x) + u_{t}(t,x) = 0,\ \ \ (t,x)\in (0,\infty)\times {\bf R}^{n},
\end{equation}
one has many previous papers, and one can cite D'Abbicco \cite{D-0}, da-Luz-Ikehata-Char\~ao \cite{Luz}, Zuazua \cite{Z}, and the references therein (For non-damped Klein-Gordon equations, one can cite several papers due to Nunes-Bastos \cite{NB} and the references therein concerning the local energy decay property and its related research).\\

By adding the mass term, the purpose of this paper is to study what kind of decay structure the equation (1.1) has. Our first result can be stated as follows.
\begin{theo}\,Let $n \geq 1$. If $[u_{0},u_{1}] \in (H^{1}({\bf R}^{n}) \cap L^{1}({\bf R}^{n})) \times (L^{2}({\bf R}^{n}) \cap L^{1}({\bf R}^{n}))$, then it is true that
\[\Vert u_{t}(t,\cdot)\Vert + \Vert u(t,\cdot)\Vert + \Vert\nabla u(t,\cdot)\Vert \leq C(1+t)^{-\frac{n}{4}}(\Vert u_{1}\Vert_{1} + \Vert u_{0}\Vert_{1}) + Ce^{-\alpha t}(\Vert u_{1}\Vert + \Vert u_{0}\Vert_{H^{1}}),\]
for any $t \geq 0$, where $C>0$ is a constant depending on $m > 0$, and $\alpha > 0$ is a independent constant of $m$.
\end{theo}
\begin{rem}{\rm If one applies the previous results introduced by \cite[Corollary 2.1 and Corollary 2.2 based on Theorem 2.2 with (i) of Hypothesis 1]{Luz}, one can easily derive
\[\Vert u_{t}(t,\cdot)\Vert + \Vert u(t,\cdot)\Vert + \Vert\nabla u(t,\cdot)\Vert \leq C t^{-\frac{n}{4}+\delta}(\Vert u_{1}\Vert_{1} + \Vert u_{0}\Vert_{1}) + Ce^{-\alpha t}(\Vert u_{1}\Vert + \Vert u_{0}\Vert_{H^{1}}),\]
for any $\delta > 0$ and $t \gg 1$. A new point of view of Theorem 1.1 is in the preciseness of the statement with $\delta = 0$ and $t \geq 0$.}
\end{rem}
Concerning the profile in asymptotic sense is read as follows.
\begin{theo}
Let $n\geq 1$. If $[u_0,u_1]\in(H^{1}({\bf R}^n) \cap L^{1,1}({\bf R}^n))\times(L^{2}({\bf R}^n)\cap L^{1,1}({\bf R}^n))$, then the solution $u(t,x)$ to problem {\rm (1.1)-(1.2)} satisfies
\[\int_{{\bf R}^n}\left\vert\mathcal{F}(u(t,\cdot))(\xi)-\left(P_1e^{-t\vert\xi\vert^2/2}\frac{\sin{(t\sqrt{\vert\xi\vert^{2} + m^{2}})}}{\sqrt{\vert\xi\vert^{2} + m^{2}}} + P_0e^{-t\vert\xi\vert^2/2}\cos{(t\sqrt{\vert\xi\vert^{2}+ m^{2}})}\right)\right\vert^2d\xi\]
\[\leq C(\Vert u_1\Vert_{1,1}^2 + \Vert u_0\Vert_{1,1}^2)(1+t)^{-\frac{n+2}{2}} + C(\Vert u_1\Vert^2+\Vert u_0\Vert_{H^{1}}^2)e^{-\omega t} + C(\Vert u_{0}\Vert_{1}^2 + \Vert u_{1}\Vert_{1}^2)e^{-\kappa t}\quad(t \geq 0),\]
where $C > 0$, $\omega > 0$, and $\kappa > 0$ are constants, which depend only on $m > 0$.
\end{theo}
\begin{rem}{\rm The obtained asymptotic profile (as $t \to \infty$) of the solution to problem (1.1)-(1.2) is still oscillate, but the effect of the oscillate part is comparatively weak as compared with the usual strongly damped wave equation case (1.3), especially in the low dimensional case for $n = 1,2$. This is because of the existence of the mass term as one can see from the factor $\frac{\sin{(t\sqrt{\vert\xi\vert^{2} + m^{2}})}}{\sqrt{\vert\xi\vert^{2} + m^{2}}}$ of the profile, that is, this factor is no effective in the low frequency region $\vert\xi\vert \ll 1$.}
\end{rem}

Based on Theorem 1.2 one can check the optimality of the decay rate just obtained in Theorem 1.1.
\begin{theo}
Let $n\geq 1$. If $[u_0,u_1]\in(H^{1}({\bf R}^n) \cap L^{1,1}({\bf R}^n))\times(L^{2}({\bf R}^n)\cap L^{1,1}({\bf R}^n))$, then the corresponding solution $u(t,x)$ to problem {\rm (1.1)-(1.2)} satisfies
\[C^{-1}(\vert P_{0}\vert + \vert P_{1}\vert)t^{-\frac{n}{4}} \leq \Vert u(t,\cdot)\Vert \leq Ct^{-\frac{n}{4}}(\Vert u_{0}\Vert_{1} + \Vert u_{1}\Vert_{1} + \Vert u_{0}\Vert_{H^{1}} + \Vert u_{1}\Vert)\]
for large $t \gg 1$, where the constant $C > 0$ depends on $m > 0$ and $n$.
\end{theo}
\begin{rem}{\rm When one compares the result of Theorem 1.3 for (1.1) with the one of (1.3), in particular, one has a significant difference in the low dimensional case $n = 1,2$, that is, in the Klein-Gordon equation case, the wave effect is extremely weak as compared with that of diffusive one, so that one can get the decay property of the $L^{2}$-norm of solutions even in the low dimensional case. This is because of the existence of the mass term $m > 0$.}
\end{rem}

In our forthcoming paper, we will study semi-linear problem to find the so called critical exponent together with dissipative structure to the equation
\begin{equation}
u_{tt}(t,x) - \Delta u(t,x)  + m^{2}u(t,x) -\Delta u_{t}(t,x) = f(u(t,x)).
\end{equation}

Our plan in this paper is as follows. In section 2, we shall prove Theorems 1.1 by the energy method in the Fourier space due to \cite{UKS}, and in section 3 we prove Theorem 1.2 by the use of the method introduced by \cite{Ik-4}. As an application, we will discuss the optimality concerning the decay rate of the $L^{2}$-norm of solutions in Section 4 to prove Theorem 1.3.\\

{\bf Notation.}\,{\small Throughout this paper, $\| \cdot\|_q$ stands for the usual $L^q({\bf R}^{n})$-norm. For simplicity of notations, in particular, we use $\| \cdot\|$ instead of $\| \cdot\|_2$. 
\[f \in L^{1,\gamma}({\bf R}^{n}) \Leftrightarrow f \in L^{1}({\bf R}^{n}), \Vert f\Vert_{1,\gamma} := \int_{{\bf R}^{n}}(1+\vert x\vert)^{\gamma}\vert f(x)\vert dx < +\infty, \quad \gamma \geq 0.\]\\
Furthermore, we denote the Fourier transform $\hat{\phi}(\xi)$ of the function $\phi(x)$ by
\begin{equation}
{\cal F}(\phi)(\xi) := \hat{\phi}(\xi) := \frac{1}{(2\pi)^{n/2}}\int_{{\bf R}^{n}}e^{-ix\cdot\xi}\phi(x)dx,
\end{equation}
where $i := \sqrt{-1}$, and $x\cdot\xi = \displaystyle{\sum_{i=1}^{n}}x_{i}\xi_{i}$ for $x = (x_{1},\cdots,x_{n})$ and $\xi = (\xi_{1},\cdots,\xi_{n})$, and the inverse Fourier transform of ${\cal F}$ is denoted by ${\cal F}^{-1}$. When we estimate several functions by applying the Fourier transform sometimes we can also use the following definition in place of (1.7)
\[{\cal F}(\phi)(\xi) := \int_{{\bf R}^{n}}e^{-ix\cdot\xi}\phi(x)dx\]
without loss of generality. We also use the notation
\[v_{t}=\frac{\partial v}{\partial t}, \quad v_{tt}=\frac{\partial^{2} v}{\partial t^{2}}, \quad \Delta = \sum^n_{i=1}\frac{\partial^2}{\partial x_i^2},\ \ x=(x_1,\cdots,x_n).\]
}


\section{Proofs of Theorems 1.1.}

To begin with, let us start with proving Theorem 1.1 by relying on the method due to \cite{UKS}.\\
In order to use the energy method in the Fourier space, we shall prepare the following notation.
\[{E_0(t,\xi):=\frac{1}{2}\vert \hat{u}_{t}\vert^2+\frac{1}{2}(\vert\xi\vert^2+m^{2})\vert\hat{u}\vert^2},\]
\[E(t,\xi):=\frac{1}{2}\vert \hat{u}_{t}\vert^2+\frac{1}{2}(\vert\xi\vert^2+m^{2})\vert\hat{u}\vert^2+\beta\rho(\xi)\mathfrak{R} (\hat{u}_t\overline{\hat{u}})+\frac{1}{2}\beta\rho(\xi)\vert\xi\vert^{2}\vert\hat{u}\vert^2\]
\[= E_{0}(t,\xi) + \beta\rho(\xi)\mathfrak{R} (\hat{u}_t\overline{\hat{u}})+\frac{1}{2}\beta\rho(\xi)\vert\xi\vert^{2\theta}\vert\hat{u}\vert^2 ,\]
\[F(t,\xi):=\vert\xi\vert^{2}\vert\hat{u}_t\vert^2+\beta\rho(\xi)(\vert\xi\vert^2+m^{2})\vert\hat{u}\vert^2,\]
\[R(t,\xi):=\beta\rho(\xi)\vert\hat{u}_t\vert^2.\]
We define a key function $\rho: {\bf R}_{\xi}^{n}\rightarrow\bf R$ by
\[\rho(\xi)=\frac{\vert\xi\vert^{2}}{1+\vert\xi\vert^{2}}.\]
The discovery of this key function is a crucial point in this section.

Now, let us apply the Fourier transform to the both sides of (1.1) together with the initial data (1.2). Then in the Fourier space ${\bf R}_{\xi}^{n}$ one has the reduced problem
\begin{equation}
\hat{u}_{tt}(t,\xi)+(\vert\xi\vert^2+m^{2})\hat{u}(t,\xi)+\vert \xi\vert^{2}\hat{u}_t(t,\xi)=0,\quad(t,\xi)\in(0,\infty)\times{\bf R}_{\xi}^{n},
\end{equation}
\begin{equation}
\hat{u}(0,\xi)=\hat{u}_0(\xi),\quad  \hat{u}_t(0,\xi)=\hat{u}_1(\xi),\quad\xi\in{\bf R}_{\xi}^{n}.
\end{equation}
Multiply both sides of (2.1) by $\overline{\hat{u}_t}$, and further $\beta \rho(\xi)\overline{\hat{u}}$. Then, by taking the real part of the resulting identities one has
\begin{equation}
\frac{d}{dt}E_0(t,\xi)+\vert\xi\vert^{2}\vert\hat{u}_t\vert^2=0,
\end{equation}
\begin{equation}
\frac{d}{dt}\{\beta\rho(\xi)\mathfrak{R}(\hat{u}_t\overline{\hat{u}})+\frac{1}{2}\beta\vert\xi\vert^{2}\rho(\xi)\vert\hat{u}\vert^2\}+\beta(\vert\xi\vert^2+m^{2})\rho(\xi)\vert\hat{u}\vert^2=\beta\rho(\xi)\vert\hat{u}_t\vert^2.
\end{equation}
By adding (2.3) and (2.4), one has
\begin{equation}
\frac{d}{dt}E(t,\xi)+F(t,\xi)=R(t,\xi).
\end{equation}
We prove
\begin{lem}  For $\beta>0$, it is true that
\[R(t,\xi)\leq\beta F(t,\xi),\quad\xi\in{\bf R}_{\xi}^{n}.\]
\end{lem}
{\it Proof.}\,Noting the facts that
\begin{equation}
\rho(\xi) \leq 1, \quad \rho(\xi) \leq \vert\xi\vert
\end{equation}
for all $\xi \in {\bf R}_{\xi}^{n}$, the statement is true.
\hfill
$\Box$
\par 
\vspace{0.2cm}
\par
It follows from (2.5) and Lemma 2.1 that
\begin{equation}
\frac{d}{dt}E(t,\xi)+(1-\beta)F(t,\xi)\leq 0,
\end{equation}
provided that the parameters $\beta>0$ are small enough such that $1-\beta > 0$.

\begin{lem}
There is a constant $M > 0$ depending only on small $\beta > 0$ such that for all $\xi\in{\bf R}_{\xi}^{n}$,
it follows that 
\[\rho(\xi)E(t,\xi)\leq M F(t,\xi).\]
\end{lem}
{\it Proof.}\,We first note the inequality
\begin{equation}
\mathfrak{R}(\hat{u}_t\overline{\hat{u}})\leq\frac{1}{2}\left(\vert\hat{u}_{t}\vert^2+\vert\hat{u}\vert^2\right).
\end{equation}
Then, it follows (2.6) and (2.8) that
\[\rho(\xi)E(t,\xi) \leq \frac{1}{2}\vert\hat{u}_{t}\vert^2 + \frac{1}{2\beta}\beta\rho(\xi)(\vert\xi\vert^{2} + m^{2})\vert\hat{u}\vert^2\]
\[+ \frac{\beta}{2}\rho(\xi)\rho(\xi)(\vert\xi\vert^{2} + m^{2})\vert\hat{u}\vert^2 + \beta\rho(\xi)\rho(\xi)\mathfrak{R}(\hat{u}_t\overline{\hat{u}})\]
\[\leq \frac{1}{2}\vert\xi\vert^{2}\vert\hat{u}_{t}\vert^2 + \frac{1}{2\beta}F(t,\xi) + \frac{1}{2}F(t,\xi) + \frac{\beta\rho(\xi)^{2}}{2}(\vert\hat{u}_{t}\vert^2+\vert\hat{u}\vert^2)\]
\[\leq \frac{1}{2}F(t,\xi) + \frac{1}{2\beta}F(t,\xi) + \frac{1}{2}F(t,\xi) + \frac{\beta}{2}\vert\xi\vert^{2}\vert\hat{u}_{t}\vert^2+ \frac{\beta}{2}\vert\xi\vert^{2}\rho(\xi)\vert\hat{u}\vert^2\]
\[\leq (1+\frac{1}{2\beta}+\frac{\beta}{2})F(t,\xi) + \frac{1}{2}\beta\rho(\xi)(\vert\xi\vert^{2} + m^{2})\vert\hat{u}\vert^2\]
\[= (1+\frac{1}{2\beta}+\frac{\beta}{2} + \frac{1}{2})F(t,\xi),\quad (\forall \xi \in {\bf R}_{\xi}^{n}),\]
which implies the desired estimate by setting
\[M :=  (\frac{1}{2\beta}+\frac{\beta}{2} + \frac{3}{2}).\]
\hfill
$\Box$
\par 
\vspace{0.2cm}
\par
Lemma 2.2 and (2.7) imply
\begin{equation}
\frac{d}{dt}E(t,\xi)+(1-\beta)\rho(\xi)M^{-1}E(t,\xi)\leq0
\end{equation}
for any $\xi\in{\bf R}_{\xi}^{n}.$ From (2.9) we find 
\begin{equation}
E(t,\xi)\leq e^{-\alpha\rho(\xi)t}E(0,\xi),\quad\xi\in {\bf R}_{\xi}^{n},
\end{equation}
where $\alpha:=(1-\beta)M^{-1} > 0$ with small $\beta \in (0,1)$.

On the other hand, in the case of $\xi\neq0$, since we have 
\begin{equation}
\pm\beta\rho(\xi)\mathfrak{R}(\hat{u}_{t}\overline{\hat{u}})\leq\frac{\beta}{2}\vert\xi\vert^{2}\rho(\xi)\vert\hat{u}\vert^2+\frac{\beta\rho(\xi)}{2}\frac{\vert\hat{u}_{t}\vert^2}{\vert\xi\vert^{2}},
\end{equation}
it follows from the definition of $E(t,\xi)$ and (2.11) with the minus sign that
\begin{equation}
E(t,\xi)\geq\frac{1}{2}\left(1-\displaystyle\frac{\beta\rho(\xi)}{\vert\xi\vert^{2}}\right)\vert\hat{u}_{t}\vert^2+\frac{1}{2}(\vert\xi\vert^2+m^{2})\vert\hat{u}\vert^2\quad (\xi \ne 0). 
\end{equation}
And also, from (2.6) again we see 
\[\frac{\beta\rho(\xi)}{\vert\xi\vert^{2}} \leq \frac{\beta}{\vert\xi\vert^{2}}\vert\xi\vert^{2} = \beta,\]
so that one has
\begin{equation}
1-\frac{\beta\rho(\xi)}{\vert\xi\vert^{2}} \geq 1-\beta >0,
\end{equation}
for $\xi\neq0$ if we choose small $\beta\in (0,1)$. So, we obtain
\begin{equation}
E(t,\xi)\geq(1-\beta)E_0(t,\xi),\quad\xi\neq0.
\end{equation}
Since $E(t,0)=E_0(t,0),$ (2.14) holds true for all $\xi\in{\bf R}_{\xi}^{n}.$ Thus, from (2.10) one has
\begin{equation}
E_0(t,\xi)\leq(1-\beta)^{-1}e^{-\alpha\rho(\xi)t}E(0,\xi),\quad\xi\in{\bf R}_{\xi}^{n}.
\end{equation}
While, because of (2.11) with plus sign, for $\xi\neq0$ one has
\[E(t,\xi)\leq\frac{1}{2}\vert\hat{u}_{t}\vert^2+\frac{1}{2}(\vert\xi\vert^2+m^{2})\vert\hat{u}\vert^2+\frac{\beta}{2}\rho(\xi)\vert\hat{u}\vert^2+\frac{\beta}{2}\rho(\xi)\vert\hat{u}_{t}\vert^2 + \frac{\beta}{2}\vert\xi\vert^{2}\vert\hat{u}\vert^2
\]
\[\leq \frac{1}{2}(1+\beta)\vert\hat{u}_{t}\vert^2 + \frac{1}{2}(\vert\xi\vert^2+m^{2})\vert\hat{u}\vert^2+ \frac{\beta}{2}(\vert\xi\vert^2+m^{2})\vert\hat{u}\vert^2 +  \frac{\beta}{2}(\vert\xi\vert^2+m^{2})\vert\hat{u}\vert^2\]
\[\leq (1+2\beta)E_{0}(t,\xi),\]
which implies 
\begin{equation}
E(t,\xi) \leq CE_{0}(t,\xi)
\end{equation}
for all $\xi \in {\bf R}_{\xi}^{n}$ with some constant $C := (1+2\beta) > 0$. By (2.15) and (2.16) with $t=0$ one has arrived at the significant estimate.
\begin{lem} Let $\beta > 0$ be a small number. Then, there is a constant $C=C(\beta)>$0 and $\alpha=\alpha(\beta)>0$ such that for all $\xi\in {\bf R}_{\xi}^{n}$ it is true that
\begin{align*}
E_0(t,\xi)\leq Ce^{-\alpha\rho(\xi)t}E_0(0,\xi).
\end{align*}
\end{lem}

{\it Proof of Theorem 1.1.}\,By lemma 2.3 and the Plancherel theorem one has
\begin{align}
\int_{{\bf R}_{x}^n}(\vert u_t\vert^2+\vert\nabla u\vert^2 + m^{2}\vert u\vert^{2})dx \leq C\int_{{\bf R}_{\xi}^n}(\vert\hat{u}_{t}\vert^2+(\vert\xi\vert^2+m^{2})\vert\hat{u}\vert^2)d\xi\leq C\int_{{\bf R}_{\xi}^n} E_0(t,\xi)d\xi\notag\\
\leq C(\int_{\vert\xi\vert\leq 1}+\int_{\vert\xi\vert\geq1})e^{-\alpha\rho(\xi)t}(\vert\hat{u}_{1}(\xi)\vert^2+(\vert\xi\vert^2+m^{2})\vert\hat{u}_{0}(\xi)\vert^2)d\xi =: C(I_{low} + I_{high}).
\end{align}
We first prepare the following standard formula.
\begin{equation}
\int_{\vert\xi\vert\leq1}e^{-\gamma\vert\xi\vert^{2}t}\vert\xi\vert^k d\xi\leq C(1+t)^{-\frac{k+n}{2}}\quad(t \geq 0),
\end{equation}
for each $k \in {\bf N}\cup\{0\}$, and $\gamma > 0$.

Now. let us start with estimating both $I_{low}$ and $I_{high}$ based on the shape of $\rho(\xi)$. In fact, in the case when $\vert\xi\vert \leq 1$, since one has
\[\rho(\xi) \geq \frac{\vert\xi\vert^{2}}{2},\]
it follows that
\[
I_{low}\leq\int_{\vert\xi\vert\leq1}e^{-\eta\vert\xi\vert^{2}t}\vert\hat{u}_{1}(\xi)\vert^2d\xi+\int_{\vert\xi\vert\leq1}e^{-\eta\vert\xi\vert^{2}t}(\vert\xi\vert^2+m^{2})\vert\hat{u}_{0}\vert^2d\xi\]
\[
\leq C\Vert u_1\Vert_{1}^{2}(1+t)^{-\frac{n}{2}}+C\Vert u_0\Vert_{1}^{2}(1+t)^{-\frac{n}{2}},
\]
where $\eta := \frac{\alpha}{2}$. Here, we have just used the formula (2.18).

Next, we shall estimate the high frequency part. First, one also has $$\vert\xi\vert\geq1\Rightarrow\frac{\vert\xi\vert^{2}}{1+\vert\xi\vert^{2}}\geq\frac{1}{2}.$$
Hence, it follows that
\[
I_{high}=\int_{\vert\xi\vert\geq1}e^{-\alpha\rho(\xi)t}(\vert\hat{u}_{1}\vert^2+(\vert\xi\vert^2+m^{2})\vert\hat{u}_{0}\vert^2)d\xi\leq\int_{\vert\xi\vert\geq1}e^{-\eta t}(\vert\hat{u}_{1}\vert^2+(\vert\xi\vert^2+m^{2})\vert\hat{u}_{0}\vert^2)d\xi
\]
\[
\leq Ce^{-\eta t}(\Vert u_{1}\Vert^{2} + \Vert u_{0}\Vert_{H^{1}}^{2}),
\]
with some constant $C > 0$. This completes the proof of Theorem 1.1.
\hfill
$\Box$
\par 
\vspace{0.2cm}
\par

\section{Proof of Theorems 1.2.}

In this section, let us prove Theorem 1.2 by employing the method due to \cite{Ik-3,Ik-4}. We first establish the asymptotic profile of the solution in the low frequency region $\vert\xi\vert \ll 1$, which is an essential ingredient. 
\begin{lem}
Let $n\geq1$. Then, it is true that there exists a constant $C>0$ depending on $m > 0$ such that for $t \geq 0$
\begin{align*}
\int_{\vert\xi\vert\leq\delta_0}&\vert\mathcal{F}(u(t,\cdot))(\xi)-\left(P_1e^{-t\vert\xi\vert^2/2}\frac{\sin{(t\sqrt{\vert\xi\vert^{2}+m^{2}})}}{\sqrt{\vert\xi\vert^{2}+m^{2}}}+P_0e^{-t\vert\xi\vert^2/2}\cos{(t\sqrt{\vert\xi\vert^{2}+m^{2}})}\right)\vert^2d\xi\\
&\leq C(\Vert u_1\Vert_{1,1}^2 + \Vert u_0\Vert_{1,1}^2)(1+t)^{-\frac{n+2}{2}}
\end{align*}
with a small positive constant $\delta_0\ll1$.
\end{lem}

In order to prove Lemma 3.1 we apply the Fourier transform with respect to the space variable $x$ of the both sides of (1.1)-(1.2). Then in the Fourier space ${\bf R}_{\xi}^n$ one has the reduced problem:
\begin{align}
\hat{u}_{tt}(t,\xi)+(\vert\xi\vert^2 + m^{2})\hat{u}(t,\xi)+\vert \xi\vert^2\hat{u}_t(t,\xi)=0,\quad(t,\xi)\in(0,\infty)\times{\bf R}_{\xi}^n,\\ 
\hat{u}(0,\xi)=\hat{u}_0(\xi),\qquad\hat{u}_t(0,\xi)=\hat {u}_1(\xi),\quad \xi\in{\bf R}_{\xi}^n.
\end{align}
Let us solve (3.1)-(3.2) directly under the condition that $0<\vert\xi\vert\leq\delta_0\ll1$. In this case we get
\begin{align}
\hat{u}(t,\xi)&=\frac{\hat{u}_1(\xi)-\sigma_2\hat{u}_0(\xi)}{\sigma_1-\sigma_2}e^{\sigma_1t}+\frac{\hat{u}_0(\xi)\sigma_1-\hat{u}_1(\xi)}{\sigma_1-\sigma_2}e^{\sigma_2t}\\
&=\frac{e^{\sigma_1t}-e^{\sigma_2t}}{\sigma_1-\sigma_2}\hat{u}_1(\xi)+\frac{\sigma_1e^{\sigma_2t}-\sigma_2e^{\sigma_1t}}{\sigma_1-\sigma_2}\hat{u}_0(\xi),
\end{align}
where $\sigma_j\in\bf {C}$ $(j=1,2)$ have forms:
\begin{align*}
\sigma_1=\frac{-\vert\xi\vert^2+i\sqrt{4(\vert\xi\vert^{2}+m^{2})-\vert\xi\vert^4}}{2}, \quad \sigma_2=\frac{-\vert\xi\vert^2-i\sqrt{4(\vert\xi\vert^{2}+m^{2})-\vert\xi\vert^4}}{2}.
\end{align*}
In this connection, the smallness $0<\vert\xi\vert\leq\delta_0\ll1$ of $\vert\xi\vert$ is assumed to guarantee 
\[D: = 4(\vert\xi\vert^{2}+m^{2}) - \vert\xi\vert^{4} > 0.\] 

Now, we use the decomposition of the initial data based on the idea due to \cite{Ik-4}:
\begin{equation}
\hat{u}_j(\xi)=A_j(\xi)-iB_j(\xi)+P_j\quad(j=0,1),
\end{equation}
where
\[A_j(\xi):=\int_{{\bf R}^n}(\cos(x\cdot\xi)-1)u_j(x)dx,\qquad B_j(x):=\int_{{\bf R}^n}\sin(x\cdot\xi)u_j(x)dx,\quad(j=0,1).\]
From  (3.4) and (3.5) we see that 
\begin{align}
\hat{u}(t,\xi)=&P_1\left(\frac{e^{\sigma_1t}-e^{\sigma_2t}}{\sigma_1-\sigma_2}\right)+P_0\left(\frac{\sigma_1e^{\sigma_2t}-\sigma_2e^{\sigma_1t}}{\sigma_1-\sigma_2}\right)\notag\\
&+(A_1(\xi)-iB_1(\xi))\left(\frac{e^{\sigma_1t}-e^{\sigma_2t}}{\sigma_1-\sigma_2}\right)\notag\\
&+(A_0(\xi)-iB_0(\xi))\left(\frac{\sigma_1e^{\sigma_2t}-\sigma_2e^{\sigma_1t}}{\sigma_1-\sigma_2}\right),
\end{align}
for all $\xi$ satisfying $0<\vert\xi\vert\leq\delta_0.$
It is easy to check that
\begin{equation}
\frac{e^{\sigma_1t}-e^{\sigma_2t}}{\sigma_1-\sigma_2}=\frac{2e^{-t\vert\xi\vert^2/2}\sin(\frac{t\sqrt{D}}{2})}{\sqrt{D}}, 
\end{equation}
\begin{equation}
\frac{\sigma_1e^{\sigma_2t}-\sigma_2e^{\sigma_1t}}{\sigma_1-\sigma_2}=\frac{\vert\xi\vert^2e^{-t\vert\xi\vert^2/2}\sin(\frac{t\sqrt{D}}{2})}{\sqrt{D}}+e^{-t\vert\xi\vert^2/2}\cos\left(\frac{t\sqrt{D}}{2}\right).
\end{equation}
If we set
\[K_1(t,\xi):=P_0\frac{\vert\xi\vert^2e^{-t\vert\xi\vert^2/2}\sin(\frac{t\sqrt{D}}{2})}{\sqrt{D}},\]
\[K_2(t,\xi):=(A_1(\xi)-iB_1(\xi))\left(\frac{e^{\sigma_1t}-e^{\sigma_2t}}{\sigma_1-\sigma_2}\right),\]
\[K_3(t,\xi):=(A_0(\xi)-iB_0(\xi))\left(\frac{\sigma_1e^{\sigma_2t}-\sigma_2e^{\sigma_1t}}{\sigma_1-\sigma_2}\right),\]
then it follows from (3.6), (3.7) and (3.8) that
\begin{align}
\hat{u}(t,\xi)=2&P_1 \frac{e^{-t\vert\xi\vert^2/2}\sin(\frac{t\sqrt{D}}{2})}{\sqrt{D}}+P_0e^{-t\vert\xi\vert^2/2}\cos\left(\frac{t\sqrt{D}}{2}\right)\notag\\
&+K_1(t,\xi)+K_2(t,\xi)+K_3(t,\xi),\quad0<\vert\xi\vert\leq\delta_0.
\end{align}
Let us shave the needless factors of the right hand side of (3.9) to get a precise shape of the asymptotic profile by using the mean value theorem. In fact, from the mean value theorem it follows that
\[\sin(\frac{t\sqrt{D}}{2}) = \sin(t\sqrt{E}) + \frac{t}{2}(\sqrt{D}-2\sqrt{E})\cos{(\epsilon(t,\xi)),}\]
\[\cos(\frac{t\sqrt{D}}{2}) = \cos(\sqrt{E}) - \frac{t}{2}(\sqrt{D}-2\sqrt{E})\sin{(\eta(t,\xi))},\]
\[\frac{2}{\sqrt{D}} = \frac{1}{\sqrt{E}} + \frac{\vert\xi\vert^{4}}{\sqrt{(4E-\vert\xi\vert^{4}\theta_{3})^{3}}}\]
with some $\theta_{3} \in (0,1)$, where 
\[E := \vert\xi\vert^{2}+ m^{2},\]
\[\epsilon(t,\xi):=\frac{t\sqrt{D}}{2}\theta_{1} + t\sqrt{E}(1-\theta_{1}),\quad\exists\theta_{1}\in(0,1),\]
\[\eta(t,\xi):= \frac{t\sqrt{D}}{2}\theta_{2} + t\sqrt{E}(1-\theta_{2})\quad\exists\theta_{2}\in(0,1),\] 
so that from (3.9) one has arrived at the identity
\[\hat{u}(t,\xi)=P_1e^{-t\vert\xi\vert^2/2}\frac{1}{\sqrt{E}}\sin(t\sqrt{E}) + P_0e^{-t\vert\xi\vert^2/2}\cos(t\sqrt{E})\]
\[+2P_1\vert\xi\vert^{4}e^{-t\vert\xi\vert^2/2}\sin(t\sqrt{E})\frac{1}{\sqrt{(4E-\vert\xi\vert^{4}\theta_{3})^{3}}} + P_{1}t e^{-t\vert\xi\vert^2/2}\frac{\sqrt{D}-2\sqrt{E}}{\sqrt{E}}\cos{(\epsilon(t,\xi))}\]
\[+P_{1}t\vert\xi\vert^{4}e^{-t\vert\xi\vert^2/2}\frac{\sqrt{D}-2\sqrt{E}}{\sqrt{(4E-\vert\xi\vert^{4}\theta_{3})^{3}}}\cos(\epsilon(t,\xi)) - \frac{P_0 t}{2}e^{-t\vert\xi\vert^2/2}(\sqrt{D}-2\sqrt{E})\sin{(\eta(t,\xi))}+\sum_{j=1}^3K_j(t,\xi).\]
\noindent
Set
\[K_4(t,\xi) := 2P_1\vert\xi\vert^{4}e^{-t\vert\xi\vert^2/2}\sin(t\sqrt{E})\frac{1}{\sqrt{(4E-\vert\xi\vert^{4}\theta_{3})^{3}}},\]
\[K_5(t,\xi) := P_{1}t e^{-t\vert\xi\vert^2/2}\frac{\sqrt{D}-2\sqrt{E}}{\sqrt{E}}\cos{(\epsilon(t,\xi))},\]
\[K_6(t,\xi) := P_{1}t\vert\xi\vert^{4}e^{-t\vert\xi\vert^2/2}\frac{\sqrt{D}-2\sqrt{E}}{\sqrt{(4E-\vert\xi\vert^{4}\theta_{3})^{3}}}\cos(\epsilon(t,\xi)),\]\\
and
\[K_{7}(t,\xi) := -\frac{P_0 t}{2}e^{-t\vert\xi\vert^2/2}(\sqrt{D}-2\sqrt{E})\sin{(\eta(t,\xi))}.\]
Then, one has arrived at the significant identity, which holds in the low frequency region $\vert\xi\vert \leq \delta_{0}$:
\begin{equation}
\hat{u}(t,\xi) = P_1 e^{-t\vert\xi\vert^2/2}\frac{\sin(t\sqrt{\vert\xi\vert^{2}+m^{2}})}{\sqrt{\vert\xi\vert^{2}+m^{2}}}+ P_0 e^{-t\vert\xi\vert^2/2}\cos(t\sqrt{\vert\xi\vert^{2}+m^{2}}) + \sum_{j=1}^7K_j(t,\xi).
\end{equation} 

In the following part, let us check decay orders of the remainder terms $K_{j}(t,\xi)$ ($j = 1,2,3,4,5,6,7$) of (3.10) based on the formula (2.18). Note that $\delta_{0} > 0$ is a sufficiently small number such that $4m^{2} > \delta_{0}^{4}$.\\
\noindent
{\it \bf(I)\,Estimate for $ K_1(t,\xi).$}
\begin{align}
\int_{\vert\xi\vert\leq\delta_0}\vert K_1(t,\xi)\vert^2d\xi&\leq\vert P_0\vert^2\int_{\vert\xi\vert\leq\delta_0}\frac{\vert\xi\vert^4 e^{-t\vert\xi\vert^2} \vert\sin{(t\sqrt{D}/2)}\vert^2}{4(\vert\xi\vert^{2}+m^{2})-\vert\xi\vert^4}d\xi\notag\\
&\leq\frac{\vert P_0\vert^2}{4m^{2}-\delta_0^4}\int_{\vert\xi\vert\leq\delta_0}\vert\xi\vert^4e^{-t\vert\xi\vert^2}d\xi\leq\frac{\vert P_0\vert^2}{4m^{2}-\delta_0^4}(1+t)^{-\frac{n+4}{2}}.
\end{align}
\noindent
{\it \bf(II)\,Estimate for $ K_4(t,\xi).$}\\
\begin{align}
\int_{\vert\xi\vert\leq\delta_0}\vert K_4(t,\xi)\vert^2d\xi&\leq 4P_{1}^{2}\int_{\vert\xi\vert\leq\delta_0}\vert\xi\vert^{8}e^{-t\vert\xi\vert^2}\vert\sin(t\sqrt{E})\vert^{2}\frac{1}{\sqrt{(4(\vert\xi\vert^{2}+m^{2}) -\vert\xi\vert^{4}\theta_{3})^{6}}}d\xi\notag\\
&\leq 4P_{1}^{2}\frac{1}{(4m^{2}-\delta_0^4)^{3}}\int_{\vert\xi\vert\leq\delta_0}\vert\xi\vert^{8}e^{-t\vert\xi\vert^2}d\xi \leq 4P_{1}^{2}\frac{1}{(4m^{2}-\delta_0^4)^{3}}(1+t)^{-\frac{n+8}{2}}.
\end{align}
Here, we shall prepare the elementary inequality such that
\begin{equation}
\vert\sqrt{D}-2\sqrt{E}\vert = \frac{\vert\xi\vert^{4}}{\sqrt{D}+2\sqrt{E}} \leq \frac{\vert\xi\vert^{4}}{\sqrt{\vert\xi\vert^{2}+m^{2}}} \leq \frac{\vert\xi\vert^{4}}{m}\quad (\vert\xi\vert \leq \delta_{0}).
\end{equation}
{\it \bf(III)\,Estimate for $ K_5(t,\xi).$}\\
It follows from (3.13) that
\begin{align}
\int_{\vert\xi\vert\leq\delta_0}\vert K_5(t,\xi)\vert^2d\xi&\leq P_{1}^{2}t^2\int_{\vert\xi\vert\leq\delta_0}e^{-t\vert\xi\vert^2}\frac{\vert\sqrt{D}-2\sqrt{E}\vert^{2}}{\vert\xi\vert^{2}+m^{2}}\vert\cos{(\varepsilon(t,\xi))}\vert^2d\xi\notag\\
&\leq m^{-4}P_{1}^{2}t^2\int_{\vert\xi\vert\leq\delta_0}e^{-t\vert\xi\vert^2}\vert\xi\vert^{8}d\xi\leq m^{-4}P_{1}^{2}(1+ t)^{-\frac{n+4}{2}}.
\end{align}
\noindent
{\it \bf(IV)\,Estimate for $ K_6(t,\xi).$}\\
It follows from (3.13) again that
\[\int_{\vert\xi\vert\leq\delta_0}\vert K_6(t,\xi)\vert^2d\xi \leq \vert P_1\vert^2 t^{2}
\int_{\vert\xi\vert\leq\delta_0}\vert\xi\vert^{8}e^{-t\vert\xi\vert^2}\frac{\vert\sqrt{D}-2\sqrt{E}\vert^{2}}{\vert 4E-\vert\xi\vert^{4}\theta_{3}\vert^{3}}\vert\cos(\varepsilon(t,\xi))\vert^{2}d\xi\]
\[\leq \frac{\vert P_1\vert^2(1+t)^{2}}{m^{2}}\vert 4m^{2}-\delta_0^4\vert^{-3}\int_{\vert\xi\vert\leq\delta_0}e^{-t\vert\xi\vert^2}\vert\xi\vert^{16}d\xi\]
\begin{equation}
 \leq \frac{\vert P_1\vert^2}{m^{2}}{\vert 4m^{2}-\delta_0^4\vert^{-3}}(1+ t)^{-\frac{n+12}{2}}.
\end{equation}
\noindent
{\it \bf(V)\,Estimate for $ K_7(t,\xi).$}\\
Because of (3.13) again one has
\begin{align}
\int_{\vert\xi\vert\leq\delta_0}\vert K_7(t,\xi)\vert^2d\xi&\leq \vert P_1\vert^2 (1+t)^{2}m^{-2}\int_{\vert\xi\vert\leq\delta_0}\vert\xi\vert^{8}e^{-t\vert\xi\vert^2}\vert\sin(\eta(t,\xi))\vert^{2}d\xi\notag\\
&\leq \vert P_1\vert^2 (1+t)^{2}m^{-2}\int_{\vert\xi\vert\leq\delta_0}\vert\xi\vert^{8}e^{-t\vert\xi\vert^2}d\xi\leq m^{-2}\vert P_1\vert^2 (1+ t)^{-\frac{n+4}{2}}.
\end{align}

In order to estimate further $K_j\hspace{1mm} (j=2,3)$, we prepare the following lemma, introduced in \cite{Ik-3}.
\begin{lem}
Let $n \geq1$. Then it holds that for all $\xi \in {\bf R}^{n}$
\[\vert A_j (\xi)\vert\leq L\vert\xi\vert\Vert u_j\Vert_{1,1}\quad(j=0,1),\]
\[\vert B_j (\xi)\vert\leq M\vert\xi\vert\Vert u_j\Vert_{1,1}\quad(j=0,1),\]
where
\[L:=\sup_{\theta\neq0}\frac{\vert1-\cos{\theta}\vert}{\vert\theta\vert}<+\infty,\quad M:=\sup_{\theta\neq0}\frac{\vert\sin{\theta}\vert}{\vert\theta\vert}<+\infty,\]
and both $A_j(\xi)$ and $B_j(\xi)$ are defined in {\rm (3.5)}.
\end{lem}
Basing on Lemma 3.2 we check decay rates of $K_{j}(t,\xi)$ with $j = 2,3$. This part is essential in this paper.\\
\noindent
{\it \bf(VI)\,Estimate for $ K_2(t,\xi).$}\\
It follows from (3.7) and Lemma 3.2 that
\begin{align}
\int_{\vert\xi\vert\leq\delta_0}\vert K_2(t,\xi)\vert^2d\xi&\leq C(L^2+M^2)\Vert u_1\Vert_{1,1}^2\int_{\vert\xi\vert\leq\delta_0}\vert\xi\vert^2 e^{-t\vert\xi\vert^2}\frac{\sin^2{(\frac{t\sqrt{D}}{2})}}{4(\vert\xi\vert^{2}+m^{2})-\vert\xi\vert^{4}}d\xi,\notag\\
&\leq C\Vert u_1\Vert_{1,1}^2(4m^{2}-\delta_{0}^{4})^{-1}\int_{\vert\xi\vert\leq\delta_0}\vert\xi\vert^{2}e^{-t\vert\xi\vert^2}d\xi\leq C\Vert u\Vert_{1,1}^2(1+t)^{-\frac{n+2}{2}}.
\end{align}
{\it \bf(VII)\,Estimate for $K_3(t,\xi)$.}\\
Similarly, 
\[\int_{\vert\xi\vert\leq\delta_0}\vert K_3(t,\xi)\vert^2d\xi \leq  C(L^2+M^2)\Vert u_0\Vert_{1,1}^2\int_{\vert\xi\vert\leq\delta_0}\vert\xi\vert^6 e^{-t\vert\xi\vert^2}\frac{\sin^2(\frac{t\sqrt{D}}{2})}{4(\vert\xi\vert^{2}+m^{2})-\vert\xi\vert^4}d\xi \]
\[+  C(L^2+M^2)\Vert u_0\Vert_{1,1}^2\int_{\vert\xi\vert\leq\delta_0}\vert\xi\vert^2 e^{-t\vert\xi\vert^2}\cos^{2}\left(\frac{t\sqrt{D}}{2}\right)d\xi\]
\[\leq C\Vert u_0\Vert_{1,1}^2\int_{\vert\xi\vert\leq\delta_0}\vert\xi\vert^6 e^{-t\vert\xi\vert^2}d\xi+C\Vert u_0\Vert_{1,1}^2\int_{\vert\xi\vert\leq\delta_0}\vert\xi\vert^2 e^{-t\vert\xi\vert^2}d\xi\]
\begin{equation}
\leq C\Vert u_0\Vert_{1,1}^2(1+t)^{-\frac{n+6}{2}}+C\Vert u_0\Vert_{1,1}^2(1+ t)^{-\frac{n+2}{2}}.
\end{equation}
Under these preparations let us prove Lemma 3.1. \\
\par
\vspace{0.2cm}
{\it Proof of Lemma 3.1.} It follows from (3.10), (3.13), (3.15), (3.16), (3.17), (3.18) and (3.18) that
\[\int_{\vert\xi\vert\leq\delta_0}\vert\mathcal{F}(u(t,\cdot))(\xi)-\left(P_1 e^{-t\vert\xi\vert^2/2}\frac{\sin(t\sqrt{\vert\xi\vert^{2}+m^{2}})}{\sqrt{\vert\xi\vert^{2}+m^{2}}}+ P_0 e^{-t\vert\xi\vert^2/2}\cos(t\sqrt{\vert\xi\vert^{2}+m^{2}}) \right)\vert^2d\xi\]
\[\leq C\Vert u_1\Vert_{1,1}^2 (1+t)^{-\frac{n+2}{2}}+C\Vert u_0\Vert_{1,1}^2(1+t)^{-\frac{n+2}{2}},\quad (t \geq 0),\]
which implies the desired estimate.
\hfill
$\Box$
\par 
\vspace{0.2cm}
\par

Final part of this section, let us prove our main Theorem 1.2 by relying on Lemma 2.3 and Lemma 3.1. Main parts are in the estimates on the high frequency region $\vert\xi\vert \gg 1$.\\
\noindent 
{\it Proof of Theorem 1.2.}
First we decompose the integrand to be estimated into two parts such that one is the low frequency part, and the other is high frequency one as follows.   
\[\int_{{\bf R}^n} \vert\mathcal{F}(u(t,\cdot))(\xi)-\left(P_{1}e^{-t\vert\xi\vert^2/2}\frac{\sin{(t\sqrt{\vert\xi\vert^{2}+m^{2}})}}{\sqrt{\vert\xi\vert^{2}+m^{2}}}+P_0e^{-t\vert\xi\vert^2/2}\cos{(t\sqrt{\vert\xi\vert^{2}+m^{2}})}\right)\vert^2d\xi\]
\[= \left(\int_{\vert\xi\vert\leq\delta_0}+\int_{\vert\xi\vert\geq\delta_0}\right)\vert \hat{u}(t,\xi)-\left(P_{1}e^{-t\vert\xi\vert^2/2}\frac{\sin{(t\sqrt{\vert\xi\vert^{2}+m^{2}})}}{\sqrt{\vert\xi\vert^{2}+m^{2}}}+P_0e^{-t\vert\xi\vert^2/2}\cos{(t\sqrt{\vert\xi\vert^{2}+m^{2}})}\right)\vert^2d\xi\]
\[=:I_l(t)+I_h(t).\]
To begin with, as a direct consequence of Lemma 3.1 we can see  
\begin{equation}
I_l(t)\leq C(\Vert u_1\Vert_{1,1}^2 + \Vert u_0\Vert_{1,1}^2)(1+t)^{-\frac{n+2}{2}}.
\end{equation}
\noindent
On the other hand, it follows from Lemma 2.3 that
\begin{equation}
\vert\hat{u}_{t}(t,\xi)\vert^2+(\vert\xi\vert^2 + m^{2})\vert\hat{u}(t,\xi)\vert^2\leq Ce^{-\alpha\rho(\xi)t}(\vert\hat{u}_{1}(\xi)\vert^2+(\vert\xi\vert^2 + m^{2})\vert\hat{u}_{0}(\xi)\vert^2),\quad\xi\in{\bf R}^n.
\end{equation}
\noindent
In this case, we see that if $\delta_0\leq\vert\xi\vert\leq 1$, then 
$$\rho(\xi) = \frac{\vert\xi\vert^{2}}{1+\vert\xi\vert^{2}}\geq\frac{\vert\xi\vert^{2}}{2}\geq\frac{\delta_0^2}{2},$$
and if $\vert\xi\vert\geq1$, then 
$$\rho(\xi) = \frac{\vert\xi\vert^{2}}{1+\vert\xi\vert^{2}}\geq\frac{\vert\xi\vert^{2}}{2\vert\xi\vert^2} = \frac{1}{2} \geq \frac{\delta_{0}^{2}}{2}.$$
Therefore,  from (3.20) it follows that
\[
m^{2}\int_{\vert\xi\vert\geq\delta_0}\vert\hat{u}(t,\xi)\vert^2d\xi \leq C\int_{\vert\xi\vert\geq\delta_0}e^{-\alpha\rho(\xi)t}\left(\vert\hat{u}_{1}(\xi)\vert^2+(\vert\xi\vert^{2}+m^{2})\vert\hat{u}_{0}(\xi)\vert^2\right)d\xi\]
\[\leq \int_{\vert\xi\vert\geq\delta_0}e^{-\alpha\frac{\delta_{0}^{2}}{2}t}\left(\vert\hat{u}_{1}(\xi)\vert^2+(\vert\xi\vert^{2}+m^{2})\vert\hat{u}_{0}(\xi)\vert^2\right)d\xi\]
\begin{equation}
\leq Ce^{-\omega t}(\Vert u_1\Vert^2+\Vert u_0\Vert_{H^{1}}^2)\quad (t \geq 0),
\end{equation}
where $\omega := \alpha\delta_{0}^{2}/2$, and the constant $C>0$ depends on $\delta_0$, so that $m > 0$.

On the other hand, 
\[
P_1^2\int_{\vert\xi\vert\geq\delta_0}e^{-t\vert\xi\vert^2}\left|\frac{\sin{(t\sqrt{\vert\xi\vert^{2}+m^{2}})}}{\sqrt{\vert\xi\vert^{2}+m^{2}}}\right|^2d\xi+P_0^2\int_{\vert\xi\vert\geq\delta_0}e^{-t\vert\xi\vert^2}\vert\cos{(t\sqrt{\vert\xi\vert^{2}+m^{2}})}\vert^2d\xi\]
\[\leq C\left(\frac{1}{m^2}+1\right)(\Vert u_1\Vert_1^2+\Vert u_0\Vert_1^2)\int_{\vert\xi\vert\geq\delta_0}e^{-t\vert\xi\vert^2}d\xi\]
\[=C\left(\frac{1}{m^2}+1\right)(\Vert u_1\Vert_1^2+\Vert u_0\Vert_1^2)\int_{\vert\xi\vert\geq\delta_0}e^{-\frac{t\vert\xi\vert^2}{2}}e^{-\frac{t\vert\xi\vert^2}{2}}d\xi\]
\[\leq C\left(\frac{1}{m^2}+1\right)(\Vert u_1\Vert_1^2+\Vert u_0\Vert_1^2)e^{-\frac{t\delta_0^2}{2}}\int_{\vert\xi\vert\geq\delta_0}e^{-\frac{t\vert\xi\vert^2}{2}}d\xi\]
\begin{equation}
\leq C(\Vert u_1\Vert_1^2+\Vert u_0\Vert_1^2)e^{-\frac{t\delta_0^2}{2}}(1+ t)^{-\frac{n}{2}}\leq C(\Vert u_1\Vert_1^2+\Vert u_0\Vert_1^2)e^{-\kappa t}\quad(t \gg 1),
\end{equation}
with some constants $C > 0$ and $\kappa > 0$. Therefore, by evaluating $I_h(t)$ based on (3.21) and (3.22), and combining it with (3.19) one can arrive at
\begin{align*}
\int_{{\bf R}^n}&\vert\mathcal{F}(u(t,\cdot))(\xi)-\{P_1e^{-t\vert\xi\vert^2/2}\frac{\sin{(t\sqrt{\vert\xi\vert^{2}+m^{2}})}}{\sqrt{\vert\xi\vert^{2}+m^{2}}}+P_0 e^{-t\vert\xi\vert^2/2}\cos{(t\sqrt{\vert\xi\vert^{2}+m^{2}})}\}\vert^2d\xi\\
&\leq C(\Vert u_1\Vert_{1,1}^2 + \Vert u_0\Vert_{1,1}^2)(1+t)^{-\frac{n+2}{2}} + Ce^{-\omega t}(\Vert u_1\Vert^2+\Vert u_0\Vert_{H^{1}}^2)\\
&\qquad+Ce^{-\kappa t}(\Vert u_{0}\Vert_{1}^2 + \Vert u_{1}\Vert_{1}^2)\quad(t \geq 0).
\end{align*}
This implies the desired estimate. 
\hfill
$\Box$
\par 
\vspace{0.2cm}
\par

\section{Optimality of decay rates.}

In this section, we shall study the optimality of the decay rates of the solution obtained in Theorem 1.1. This is studied based on the results obtained in Theorems 1.1 and 1.2.\\
We first compute two simple estimates from above on the quantities below.\\
Set 
\[I_{1}(t) := \int_{{\bf R}_{\xi}^{n}}e^{-t\vert\xi\vert^2}\left\vert \frac{\sin{(t\sqrt{\vert\xi\vert^{2}+m^{2}})}}{\sqrt{\vert\xi\vert^{2}+m^{2}}}\right\vert^{2}d\xi,\]
\[I_{0}(t) := \int_{{\bf R}_{\xi}^{n}}e^{-t\vert\xi\vert^2}\vert\cos(t\sqrt{\vert\xi\vert^{2}+m^{2}})\vert^2d\xi.\]
Now, it follows that
\begin{equation}
I_{1}(t) \leq \frac{1}{m^{2}}\int_{{\bf R}_{\xi}^{n}}e^{-t\vert\xi\vert^2}d\xi \leq Cm^{-2}(1+t)^{-\frac{n}{2}},
\end{equation}
and
\begin{equation}
I_{0}(t) \leq \int_{{\bf R}_{\xi}^{n}}e^{-t\vert\xi\vert^2}d\xi \leq C(1+t)^{-\frac{n}{2}}.
\end{equation}
Next, in order to get the estimates from below about $I_{j}(t)$ ($j = 0,1$), we use the polar co-ordinate transformation as follows:
\[I_{1}(t) = \vert\omega_{n}\vert t^{-\frac{n}{2}}\int_{0}^{\infty}\sigma^{n-1}e^{-\sigma^{2}}\frac{t\sin^{2}\left(t\sqrt{\frac{\sigma^{2}}{t}+m^{2}}\right)}{tm^{2} + \sigma^{2}}d\sigma\]
\begin{equation}
= \frac{\vert\omega_{n}\vert}{2}t^{-\frac{n}{2}}\int_{0}^{\infty}\sigma^{n-1}e^{-\sigma^{2}}\frac{t}{\sigma^{2}+tm^{2}}d\sigma - \frac{\vert\omega_{n}\vert}{2}t^{-\frac{n}{2}}\int_{0}^{\infty}\sigma^{n-1}e^{-\sigma^{2}}\frac{t}{\sigma^{2}+tm^{2}}\cos(2t\sqrt{\frac{\sigma^{2}}{t} + m^{2}})d\sigma,
\end{equation}
where the surface measure of the $n$-dimensional unit ball is defined by
\[\vert\omega_{n}\vert := \int_{\vert\omega\vert = 1}d\omega.\]
Let us estimate (4.3) as follows. Since
\[0 \leq \sigma^{n-1}e^{-\sigma^{2}}\frac{t}{\sigma^{2}+tm^{2}} \leq \frac{1}{m^{2}}\sigma^{n-1}e^{-\sigma^{2}} \in L^{1}(0,\infty),\]
and
\[\lim_{t \to \infty}\sigma^{n-1}e^{-\sigma^{2}}\frac{t}{\sigma^{2}+tm^{2}} = \frac{1}{m^{2}}\sigma^{n-1}e^{-\sigma^{2}},\]
it follows from the Lebesgue convergence theorem that
\[\lim_{t \to \infty}\int_{0}^{\infty}\sigma^{n-1}e^{-\sigma^{2}}\frac{t}{\sigma^{2}+tm^{2}}d\sigma = \frac{1}{m^{2}}\int_{0}^{\infty}\sigma^{n-1}e^{-\sigma^{2}}d\sigma =: K_{n} > 0,\]
so that there exists a constant $t_{0} \gg 1$ such that for all $t \geq t_{0}$ one has
\[\int_{0}^{\infty}\sigma^{n-1}e^{-\sigma^{2}}\frac{t}{\sigma^{2}+tm^{2}}d\sigma \geq \frac{1}{2}K_{n}.\]
This implies
\begin{equation}
I_{1}(t) \geq \frac{K_{n}\vert\omega_{n}\vert}{4}t^{-\frac{n}{2}} - \frac{\vert\omega_{n}\vert}{2}t^{-\frac{n}{2}}\int_{0}^{\infty}\sigma^{n-1}e^{-\sigma^{2}}\frac{t}{\sigma^{2}+tm^{2}}\cos(2t\sqrt{\frac{\sigma^{2}}{t} + m^{2}})d\sigma
\end{equation}
for all $t \geq t_{0}$. Set
\[C_{1}(t) := \int_{0}^{\infty}\sigma^{n-1}e^{-\sigma^{2}}\frac{t}{\sigma^{2}+tm^{2}}\cos(2t\sqrt{\frac{\sigma^{2}}{t} + m^{2}})d\sigma.\]
Then, it holds that
\[C_{1}(t)\]
\begin{equation}
= \frac{1}{m^{2}}\int_{0}^{\infty}\sigma^{n-1}e^{-\sigma^{2}}\cos(2t\sqrt{\frac{\sigma^{2}}{t} + m^{2}})d\sigma - \frac{1}{m^{2}}\int_{0}^{\infty}\sigma^{n-1}e^{-\sigma^{2}}\frac{\sigma^{2}}{\sigma^{2}+tm^{2}}\cos(2t\sqrt{\frac{\sigma^{2}}{t} + m^{2}})d\sigma.
\end{equation}
Here, by the Lebesgue convergence theorem it follows again that
\begin{equation}
0 \leq \int_{0}^{\infty}\sigma^{n-1}e^{-\sigma^{2}}\frac{\sigma^{2}}{\sigma^{2}+t m^{2}}\vert \cos(2t\sqrt{\frac{\sigma^{2}}{t} + m^{2}})\vert d\sigma \leq \int_{0}^{\infty}\sigma^{n-1}e^{-\sigma^{2}}\frac{\sigma^{2}}{\sigma^{2}+t m^{2}}d\sigma \to 0
\end{equation}
as $t \to \infty$. Furthermore, it follows from the Lebesgue theorem that (see Appendix below)
\begin{equation}
\lim_{t \to \infty}\int_{0}^{\infty}\sigma^{n-1}e^{-\sigma^{2}}\cos(2t\sqrt{\frac{\sigma^{2}}{t} + m^{2}})d\sigma = 0.
\end{equation} 
Therefore, it follows from (4.5)-(4.7) that
\begin{equation}
\lim_{t \to \infty}C_{1}(t) = 0.
\end{equation}
Because of (4.4) and (4.8) one has
\[I_{1}(t) \geq \frac{K_{n}\vert\omega_{n}\vert}{4}t^{-\frac{n}{2}} - \frac{\vert\omega_{n}\vert}{2}t^{-\frac{n}{2}}o(1)\quad (t \to \infty).\]
Thus, there exists a constant $t_{1} (\geq t_{0})$ such that for all $t \geq t_{1}$
\begin{equation}
I_{1}(t) \geq \frac{K_{n}\vert\omega_{n}\vert}{8}t^{-\frac{n}{2}}.
\end{equation}

Next, concerning $I_{0}(t)$, similarly to the computation for $I_{1}(t)$ one can get
\[I_{0}(t) = \frac{\vert\omega_{n}\vert}{2}t^{-\frac{n}{2}}\int_{0}^{\infty}\sigma^{n-1}e^{-\sigma^{2}}d\sigma + \frac{\vert\omega_{n}\vert}{2}t^{-\frac{n}{2}}\int_{0}^{\infty}\sigma^{n-1}e^{-\sigma^{2}}\cos(2t\sqrt{\frac{\sigma^{2}}{t} + m^{2}})d\sigma,\]
so that by the Riemann-Lebesgue theorem it follows that
\begin{equation}
I_{0}(t) \geq \frac{L_{n}\vert\omega_{n}\vert}{4}t^{-\frac{n}{2}},\quad (t \geq t_{1})
\end{equation}
where
\[L_{n} := \int_{0}^{\infty}\sigma^{n-1}e^{-\sigma^{2}}d\sigma.\]
\noindent
On the other hand, set
\[I_{2}(t) := \int_{{\bf R}_{\xi}^{n}}\frac{e^{-t\vert\xi\vert^{2}}}{\sqrt{\vert\xi\vert^{2}+m^{2}}}\sin\left(2t\sqrt{\vert\xi\vert^{2}+m^{2}}\right)d\xi.\]
Then, by the polar-coordinate transform one has
\begin{equation}
I_{2}(t) = \vert\omega_{n}\vert t^{-\frac{n}{2}}\int_{0}^{\infty}\frac{e^{-\sigma^{2}}\sigma^{n-1}}{\sqrt{(\sigma^{2}/t) + m^{2}}}\sin(2t\sqrt{(\sigma^{2}/t) + m^{2}}) d\sigma.
\end{equation}
Since $e^{-\sigma^{2}}\sigma^{k} \in L^{1}(0,\infty)$ for $k \in {\bf N}\cup\{0\}$, one can also get 
\[\int_{0}^{\infty}\frac{e^{-\sigma^{2}}\sigma^{n-1}}{\sqrt{(\sigma^{2}/t) + m^{2}}}\sin(2t\sqrt{(\sigma^{2}/t) + m^{2}}) d\sigma = (o(1) + O(t^{-1}))\quad (t \to \infty),\]
so that from (4.11) one can get
\begin{equation}
I_{2}(t) = \vert\omega_{n}\vert t^{-\frac{n}{2}}(o(1) + O(t^{-1}))\quad (t \to \infty).
\end{equation}
Indeed, because of the mean value theorem it follows that
\[\frac{1}{\sqrt{\frac{\sigma^{2}}{t}+m^{2}}} = \frac{1}{m}-\frac{\sigma^{2}}{2t}\frac{1}{((\theta\sigma^{2}/t)+m^{2})^{3/2}}\]
for some $\theta \in (0,1)$, so that one has 
\[\int_{0}^{\infty}\frac{e^{-\sigma^{2}}\sigma^{n-1}}{\sqrt{(\sigma^{2}/t) + m^{2}}}\sin(2t\sqrt{(\sigma^{2}/t) + m^{2}}) d\sigma\]
\[= \frac{1}{m}\int_{0}^{\infty}e^{-\sigma^{2}}\sigma^{n-1}\sin(2t\sqrt{(\sigma^{2}/t) + m^{2}}) d\sigma\]
\[-\frac{1}{2t}\int_{0}^{\infty}e^{-\sigma^{2}}\sigma^{n+1}\frac{1}{((\theta\sigma^{2}/t) +m^{2})^{3/2}}\sin(2t\sqrt{(\sigma^{2}/t) + m^{2}}) d\sigma.\]
Here, one can estimate as follows:
\[\frac{1}{2t}\int_{0}^{\infty}e^{-\sigma^{2}}\sigma^{n+1}\frac{1}{((\theta\sigma^{2}/t) +m^{2})^{3/2}}\vert\sin(2t\sqrt{(\sigma^{2}/t) + m^{2}})\vert d\sigma\]
\[\leq \frac{1}{2m^{3}t}\int_{0}^{\infty}e^{-\sigma^{2}}\sigma^{n+1}d\sigma = O(t^{-1})\quad (t \to \infty),\]
so that 
\[I_{2}(t) = \vert\omega_{n}\vert t^{-\frac{n}{2}}\left(\frac{1}{m}\int_{0}^{\infty}e^{-\sigma^{2}}\sigma^{n-1}\sin(2t\sqrt{(\sigma^{2}/t) + m^{2}}) d\sigma + O(t^{-1}) \right),\]
which implies the validity of (4.12), where one has just used the similar argument to deriving (4.7) (see also Appendix below).

Now, noting the identity
\[\Vert P_{1}e^{-\frac{t\vert\xi\vert^{2}}{2}}\frac{\sin(t\sqrt{\vert\xi\vert^{2}+m^{2}})}{\sqrt{\vert\xi\vert^{2}+m^{2}}}+ P_{0}e^{-\frac{t\vert\xi\vert^{2}}{2}}\cos(t\sqrt{\vert\xi\vert^{2}+m^{2}}) \Vert^{2}\]
\[= \vert P_{1}\vert^{2}I_{1}(t) + \vert P_{0}\vert^{2}I_{0}(t) + P_{1}P_{0}I_{2}(t),\]
it follows from (4.9), (4.10) and (4.12) that
\[\Vert P_{1}e^{-\frac{t\vert\xi\vert^{2}}{2}}\frac{\sin(t\sqrt{\vert\xi\vert^{2}+m^{2}})}{\sqrt{\vert\xi\vert^{2}+m^{2}}}+ P_{0}e^{-\frac{t\vert\xi\vert^{2}}{2}}\cos(t\sqrt{\vert\xi\vert^{2}+m^{2}}) \Vert^{2}\]
\[\geq \frac{K_{n}\vert\omega_{n}\vert}{8}\vert P_{1}\vert^{2}t^{-\frac{n}{2}} +\frac{L_{n}\vert\omega_{n}\vert}{4}\vert P_{0}\vert^{2}t^{-\frac{n}{2}} + P_{1}P_{0}\vert\omega_{n}\vert t^{-\frac{n}{2}}(o(1) + O(t^{-1}))\]
for large $t \gg 1$. This shows
\[\Vert P_{1}e^{-\frac{t\vert\xi\vert^{2}}{2}}\frac{\sin(t\sqrt{\vert\xi\vert^{2}+m^{2}})}{\sqrt{\vert\xi\vert^{2}+m^{2}}}+ P_{0}e^{-\frac{t\vert\xi\vert^{2}}{2}}\cos(t\sqrt{\vert\xi\vert^{2}+m^{2}}) \Vert^{2}\]
\begin{equation}
\geq \vert P_{1}\vert^{2}\frac{K_{n}\vert\omega_{n}\vert}{16}t^{-\frac{n}{2}} +\vert P_{0}\vert^{2}\frac{L_{n}\vert\omega_{n}\vert}{8}t^{-\frac{n}{2}}
\end{equation}
for large $t \gg 1$.\\

Based on the estimates (4.13) and Theorems 1.1 and 1.2 combined with the Plancherel theorem, one can prove Theorem 1.3 as follows. This is rather standard nowadays. Indeed, one has\\
\[\Vert u(t,\cdot)\Vert = \Vert\hat{u}(t,\cdot)\Vert \geq \Vert P_{1}e^{-\frac{t\vert\xi\vert^{2}}{2}}\frac{\sin(t\sqrt{\vert\xi\vert^{2}+m^{2}})}{\sqrt{\vert\xi\vert^{2}+m^{2}}}+ P_{0}e^{-\frac{t\vert\xi\vert^{2}}{2}}\cos(t\sqrt{\vert\xi\vert^{2}+m^{2}})\Vert\]
\[-\Vert \hat{u}(t,\cdot) - \left(P_{1}e^{-\frac{t\vert\xi\vert^{2}}{2}}\frac{\sin(t\sqrt{\vert\xi\vert^{2}+m^{2}})}{\sqrt{\vert\xi\vert^{2}+m^{2}}}+ P_{0}e^{-\frac{t\vert\xi\vert^{2}}{2}}\cos(t\sqrt{\vert\xi\vert^{2}+m^{2}}) \right)\Vert\]
\[\geq  (\vert P_{1}\vert^{2}\frac{K_{n}\vert\omega_{n}\vert}{16} + \vert P_{0}\vert^{2}\frac{L_{n}\vert\omega_{n}\vert}{8})^{1/2}t^{-\frac{n}{4}} + O(t^{-\frac{n+2}{4}})\]
for large $t \gg 1$. Anyway, under the assumption of Theorem 1.3, one can get the estimate of $\Vert u(t,\cdot)\Vert$ from below. 
\noindent
Concerning the estimates from above, one can use Theorem 1.1 directly. This completes the proof of Theorem 1.3.

\hfill
$\Box$
\par

\section{Appendix.}

In this section, we will check the validity of (4.7) by using the Lebesgue convergence theorem. In this connection, in the case when $m = 0$ (massless type) (4.7) is a direct consequence of the Riemann-Lebesgue theorem because of $\sigma^{n-1}e^{-\sigma^{2}} \in L^{1}(0,\infty)$. In the case of $m > 0$, one needs a little troublesome computations by using a specialty of the functions $\sigma^{n-1}e^{-\sigma^{2}}$ as follows.\\
We set
\[
I(t) := \int_{0}^{\infty}\sigma^{n-1}e^{-\sigma^{2}}\cos(2t\sqrt{\frac{\sigma^{2}}{t} + m^{2}})d\sigma.
\]
Then, one has
\[I(t) = t^{\frac{n}{2}}\int_{0}^{\infty}\sigma^{n-1}e^{-t\sigma^{2}}\cos\left(2t\sqrt{\sigma^{2}+m^{2}}\right)d\sigma,\]
so that one can get
\[
0 \leq I(t) \leq t^{\frac{n}{2}}\int_{0}^{\infty}\sigma^{n-1}e^{-t\sigma^{2}}d\sigma = \left(\int_{0}^{1}+\int_{1}^{\infty} \right) t^{\frac{n}{2}}\sigma^{n-1}e^{-t\sigma^{2}}d\sigma =: I_{low}(t) + I_{high}(t).
\]
Let us estimate $ I_{low}(t)$ and $I_{high}(t)$, respectively.\\

To begin with, one can get
\begin{equation}\lim_{t \to \infty} t^{\frac{n}{2}}\sigma^{n-1}e^{-t\sigma^{2}} = 0\quad (a.e. \sigma \in [0,\infty))\end{equation}
So, in order to apply the Lebesgue dominated convergence theorem, one has to get a priori bounds for the $t$-dependent functions $t^{\frac{n}{2}}\sigma^{n-1}e^{-t\sigma^{2}}$ in each interval $[0,1]$ and $[1,\infty)$.\\
\noindent 
{\bf Case I}: In the case of $\sigma \in (0,1]$, there exists a constant $M > 0$ such that
\[0 \leq t^{\frac{n}{2}}\sigma^{n-1}e^{-t\sigma^{2}} \leq M,\]
where $M$ is independent of any $t \gg 1$ and $\sigma \in (0,1]$. Together with (5.1), one can get
\[\lim_{t \to \infty}I_{low}(t) = 0.\]
\\
\noindent 
{\bf Case II}: In the case when $\sigma \in [1,\infty)$, if we choose $\ell \in {\bf N}$ large enough such that $2\ell > n$, then one can get
\[0 \leq t^{\frac{n}{2}}\sigma^{n-1}e^{-t\sigma^{2}} \leq \frac{\sigma^{n-1}t^{\frac{n}{2}}\ell !}{\ell ! + (t\sigma^{2})^{\ell}} \leq t^{-(\ell-\frac{n}{2})}\ell !\sigma^{-2\ell + n-1} \leq \ell !\sigma^{-2\ell+n-1}\]
with $\sigma^{-2\ell + n-1} \in L^{1}(1,\infty)$ for $t \gg 1$.Together with (5.1), one can also get
\[\lim_{t \to \infty}I_{high}(t) = 0.\]\\
\noindent
Anyway, from the argument above, one can check the validity of (4.7).

\hfill
$\Box$
\par 
\vspace{0.2cm}
\par

\noindent{\em Acknowledgement.}
\smallskip
The work of the author (R. IKEHATA) was supported in part by Grant-in-Aid for Scientific Research (C) 15K04958 of JSPS.



\begin{thebibliography}{99}
\bibitem{CLI} R. C. Char\~ao, C. R. da Luz and R. Ikehata, Sharp decay rates for wave equations with a fractional damping via new method in the Fourier space, J. Math Anal. Appl. 408 (1)(2013), 247-255.
\bibitem{D-0} M. D'Abbicco, Asymptotics for damped evolution operators with mass-like terms, Complex Analysis and Dynamical Systems VI, Contemporary Mathematics 653 (2015), 93-116.
\bibitem{DE} M. D'Abbicco and M. R. Ebert, Diffusion phenomena for the wave equation with structural damping in the $L^{p}$-$L^{q}$ framework, J. Diff. Eqns 256 (7) (2014), 2307-2336.
\bibitem{DEP} M. D'Abbicco, M. R. Ebert and T. Picon, Long time decay estimates in real Hardy spaces for evolution equations with structural dissipation, J. Pseudo-Differ. Oper. Appl. 7 (2016), 261-293. ISSN: 1662-9981; http://dx.doi.org/10.1007/s11868-015-0141-9
\bibitem{DR} M. D'Abbicco and M. Reissig, Semilinear structural damped waves, Math. Methods Appl. Sci. 37(11)(2014), 1570-1592.
\bibitem{Luz} C. R. da Luz, R. Ikehata and R. C. Char\~ao, Asymptotic behavior for abstract evolution differential equations of second order, J. Diff. Eqns 259 (2015), 5017-5039.
\bibitem{HO} T. Hosono and T. Ogawa, Large time behavior and  $L^{p}$-$L^{q}$ estimate of $2$-dimensional nonlinear damped wave equations, J. Diff. Eqns 203 (2004), 82-118.
\bibitem{Ik-3} R. Ikehata, New decay estimates for linear damped wave equations and its application to nonlinear problem, Math. Meth. Appl. Sci. 27 (2004), 865-889.
\bibitem{Ik-4} R. Ikehata, Asymptotic profiles for wave equations with strong damping, J. Diff. Eqns 257 (2014), 2159-2177.
\bibitem{INatsume} R. Ikehata and M. Natsume, Energy decay estimates for wave equations with a fractional  damping, Diff. Int. Eqns 25 (9-10) (2012), 939-956.
\bibitem{IO} R. Ikehata and M. Onodera, Remarks on large time behavior of the $L^{2}$-norm of solutions to strongly damped wave equations, Diff. Int. Eqns 30 (2017), 505-520.
\bibitem{IT} R. Ikehata and H. Takeda, Asymptotic profiles of solutions for structural damped wave equations, arXiv: 1607.01839v1 [math. AP] 6 Jul 2016.
\bibitem{ITY} R. Ikehata, G. Todorova and B. Yordanov, Wave equations with strong damping in Hilbert spaces, J. Diff. Eqns 254 (2013), 3352-3368.
\bibitem{K} G. Karch, Selfsimilar profiles in large time asymptotics of solutions to damped wave equations, Studia Math. 143 (2000), 175-197.
\bibitem{M} H. Michihisa, Expanding methods for evolution operators of strongly damped wave equations, preprint (2018).
\bibitem{Na} T. Narazaki, $L^{p}$-$L^{q}$ estimates for damped wave equations and their applications to semilinear problem, J. Math. Soc. Japan 56 (2004), 585-626.
\bibitem{N-2} K. Nishihara, $L^{p}$-$L^{q}$ estimates to the damped wave equation in $3$-dimensional space and their application, Math. Z. 244 (2003), 631-649.
\bibitem{NB} R. S. O. Nunes and W. D. Bastos, Energy decay for the linear Klein-Gordon equation and boundary control, J. Math. Anal. Appl. 414 (2014), 934-944.
\bibitem{p} G. Ponce, Global existence of small solutions to a class of nonlinear evolution equations, Nonlinear Anal. 9(5) (19), 399-418.
\bibitem{Shibata} Y. Shibata, On the rate of decay of solutions to linear viscoelastic equation, Math. Meth. Appl. Sci. 23 (2000), 203-226.
\bibitem{T} H. Takeda, Higher-order expansion of solutions for a damped wave equation, Asymptotic Anal. 94 (2015), 1-31. DOI: 10.3233/ASY-151295
\bibitem{UKS} T. Umeda, S. Kawashima and Y. Shizuta, On the decay of solutions to the linearized equations of electro-magneto-fluid dynamics, Japan J. Appl. Math. 1 (1984), 435-457.
\bibitem{W} Y. Wakasugi, On diffusion phenomena for the linear wave equation with space-dependent damping, J. Hyperbolic Diff. Eqns 11 (2014), 795-819.
\bibitem{Z} E. Zuazua, Exponential decay for the semilinear wave equation with localized damping in unbounded domains, J. Math. Pures et Appl. 70 (1992), 513-529.

\end{thebibliography}
\end{document}